\def\ez{E^{\zeta}}
\def\ex{E^\xi}
\def\cza{C^{\zeta}_a}
\def\czb{C^{\zeta}_b}
\def\z{\zeta}
\def\l{\lambda}
\def\fzab{f^\z_{a,b}}
\def\fxab{f^\xi_{a,b}}
\def\av{\vec{a}}
\def\bv{\vec{b}}
\def\max{{\rm max}}
\def\min{{\rm min}}
\def\uat{u^a_t}
\def\utb{u^b_t}
\def\vat{v^a_t}
\def\vat{v^a_t}\def\qed{$\Box$\medskip}
\def\tzi{\theta^{\zeta}_i}
\def\tei{\theta^{\epsilon}_i}
\def\tzii{\bar{\theta}^{\zeta}_i}
\def\tgi{\theta^{\gamma}_i}
\newtheorem{theorem}{Theorem}
\newtheorem{lemma}[theorem]{Lemma}

\newtheorem{proposition}[theorem]{Proposition}
\newtheorem{corollary}[theorem]{Corollary}
\newtheorem{definition}[theorem]{Definition}
\def\proof{\noindent{\bf Proof.}\hspace{2mm}}
\documentclass[12pt]{article}
\usepackage{amsfonts}
\usepackage{amssymb}
\usepackage{latexsym}

\newif\ifdoublespace
\doublespacetrue

\ifdoublespace

\newcommand{\bibsp}{\vspace*{.6\baselineskip}}
\else

\newcommand{\bibsp}{\relax}
\fi
\def\gap{\langle \aleph_0,\aleph_1\rangle \rightarrow \langle 
\lambda,\lambda^+ \rangle}
\author{Juliette Kennedy
\thanks{Research partially supported by 
grant 1011049 of the Academy of Finland.}\\
Department of Mathematics\\
University of Helsinki\\
Helsinki, Finland
\and
Saharon Shelah
\thanks{Research
 partially supported by the Binational Science Foundation.
 Publication number 769.}\\
Institute of Mathematics\\
Hebrew University\\
Jerusalem, Israel\\}
\title{On 
regular reduced products\thanks{This paper was written while the
authors were guests of the Mittag-Leffler Institute,
Djursholm, Sweden. The authors are grateful to the
Institute for its support.}}

\begin{document}

\maketitle

\begin{abstract}
Assume \(\gap\).
Assume $M$ is a model of a first order  theory $T$ of
cardinality at 
most $\lambda^+$ in a vocabulary $\mathcal{L}(T)$
of cardinality $\leq
\lambda$. Let $N$ be a model with the same vocabulary. 
Let $\Delta$ be a set of 
first order formulas in $\mathcal{L}(T)$ and let $D$ be a 
regular filter on 
$\lambda$. Then 
 $M$ is $\Delta$-embeddable into the reduced power 
$N^{\lambda}/ D$, provided that every 
$\Delta$-existential formula true in $M$ is true also in $N$.
We obtain the following corollary:
for $M$ as above and $D$ a regular ultrafilter over $\lambda$, 
$M^{\lambda}/ D$ is $\lambda ^{++}$-universal.
Our second result is as follows:
For $i<\mu$ let $M_i$ and $N_i$ be elementarily
equivalent models of a vocabulary which has
has cardinality $\le\l$. Suppose $D$ is a regular
filter on $\mu$ and \(\gap\) holds. 
We show that then the second player has
a winning strategy in the Ehrenfeucht-Fraisse
game of length $\l^+$ on $\prod_iM_i/D$
and $\prod_iN_i/D$. This yields the following corollary:
Assume GCH and $\l$ regular (or just $\gap$
and $2^\l=\l^+$). 
For $L$, $M_i$ and $N_i$ be
as above, if $D$ is a regular filter on $\l$,
then  $\prod_iM_i/D\cong\prod_iN_i/D$.
\end{abstract}

\section{Introduction}

Suppose \(M\) is a first order structure
and \(F\) is the Frechet filter on \(\omega\).
Then the reduced power \(M^{\omega}/F\)
is \(\aleph_1\)-saturated and hence 
\(\aleph_2\)-universal (\cite{jo}).
This was generalized by Shelah in \(\cite{shb}\)
to any filter \(F\) on
\(\omega\) for which \(B^\omega/F\) is
\(\aleph_1\)-saturated, where \(B\) is the
two element Boolean algebra,
and in \cite{ks} to all regular filters on 
$\omega$.
In the first part of this paper we use the combinatorial principle
\(\Box^{b^*}_\lambda\) of Shelah~\cite{Sh269} to
generalize the result
from \(\omega\) to arbitrary \(\lambda\), assuming $\gap$.
This gives a partial solution to Conjecture 19
in \cite{CK}: if $D$ is a regular ultrafilter over $\lambda$,
then for all infinite $M$, the ultrapower $M^{\lambda}/D$ is
$\lambda^{++}$-universal. 

The second part of this paper addresses Problem 18
in \cite{CK}, which  asks if it is true that if 
$D$ is a regular ultrafilter 
over $\lambda$,
then for all elementarily equivalent 
models $M$ and $N$ of cardinality $\le\l$
in a vocabulary of cardinality $\le\l$, the 
ultrapowers $M^{\lambda}/D$ 
and $N^{\lambda}/D$ are isomorphic.
Keisler \cite{k} proved this for good $D$ 
assuming $2^\l=\l^+$. Benda \cite{b} weakened
"good" to "contains a good filter".
We prove the claim in full generality,
assuming $2^\l=\l^+$ and $\gap$.

Regarding our assumption $\gap$,
by Chang's Two-Cardinal Theorem (\cite{chang})
$\gap$ is a consequence of $\l=\l^{<\l}$. So our 
Theorem~\ref{first} settles
Conjecture 19 of \cite{CK},
and Theorem~\ref{second} settles
Conjecture 18 of \cite{CK},  under GCH for $\l$ regular. For singular
strong limit cardinals $\gap$ follows from $\Box_\l$ (Jensen
\cite{Jen}). In the so-called Mitchell's model (\cite{Mit})
$\langle \aleph_0,\aleph_1\rangle\not\rightarrow
\langle\aleph_1,\aleph_2\rangle$, so our assumption
is independent of ZFC.

\section{Universality}

\begin{definition}
Suppose $\Delta$ is a set of first order formulas of
vocabulary $L$. The set
of {\em $\Delta$-existential} formulas is the set of
 formulas of the form
$$\exists x_1...\exists x_n(\phi_1\wedge...\wedge\phi_n),$$
where each $\phi_i$ is in $\Delta$. The set
of {\em weakly $\Delta$-existential} formulas is the set of
 formulas of the above form,
where each $\phi_i$ is in $\Delta$ or is the negation of a
formula in $\Delta$.
If $M$ and $N$ are $L$-structures
and $h:M\rightarrow N$, we say that $h$ is a 
{\em $\Delta$-homomorphism} if $h$ preserves the
truth of $\Delta$-formulas. If $h$ preserves also
the truth of negations of $\Delta$-formulas,
it is called a {\em $\Delta$-embedding}.
\end{definition}

\begin{theorem} \label{first} Assume $\gap$.
Let $M$ be a model of a first order  theory $T$ of
cardinality at 
most $\lambda^+$, in a language $L$ of cardinality $\leq
\lambda$ and let \(N\) be a model with the
same vocabulary. Let $\Delta$ be a set of 
first order formulas in $L$ and let $D$ be a regular 
filter on $\lambda$. We assume that every weakly 
$\Delta$-existential sentence true in $M$ is true also in $N$.
Then there 
is a $\Delta$-embedding of $M$ into the reduced power 
$N^{\lambda}/ D$. 
\end{theorem}

By letting $\Delta$ be the set of all first order sentences,
we get from Theorem~\ref{first} and \L o\'s' Lemma:
 
\begin{corollary} Assume $\gap$.
If $M$ is a model with vocabulary $\leq\lambda$, and $D$
is a regular ultrafilter over $\lambda$, then $M^{\lambda}/D$
is $\lambda^{++}$-universal, i.e. if $M'$ is of cardinality 
$\leq \lambda^+$, and $M'\equiv M$, then $M'$ is elementarily
embeddable into the ultrapower $M^{\lambda}/D$. 
\end{corollary}

We can replace "weakly $\Delta$-existential" by
"$\Delta$-existential" in the Theorem, if we only want
a $\Delta$-homomorphism.

The proof of Theorem \ref{first} is an
induction over $\lambda$ and $\lambda^+$ respectively, as follows.
Suppose $M=\{a_\zeta:\zeta<\lambda^+\}$.
We associate to  each \(\zeta<\lambda^+\) 
finite sets $u^{\zeta}_i$, 
$i<\lambda$, and represent the
formula set $\Delta$ as a union of finite
sets $\Delta_i$. At stage $i$, for each $\zeta <\lambda^+$ we 
consider the $\Delta_i$-type 
of the elements $a_{\zeta}$ of the model whose indices lie in 
the set $u^{\zeta}_i$, $\zeta <\lambda^+$. 
This will yield a witness $f_{\zeta}(i)$ in $N$ at stage $i,\zeta$.
Our embedding is then given by 
$a_{\zeta} \mapsto \langle f_{\zeta}(i):
i<\lambda \rangle /D$.

We need first an important lemma, reminiscent of Proposition
5.1 in \cite{Sh269}:

\def\uai{u^{\alpha}_i}
\def\uzi{u^{\zeta}_i}
\def\uei{u^{\epsilon}_i}
\def\ugi{u^{\gamma}_i}

\begin{lemma}\label{main}
Assume \(\gap\). Let $D$ be a regular filter on $\l$.
There exist sets \(\uzi\) and integers $n_i$
for each $\zeta<\lambda^+$ and $i<\lambda$
such that for each \(i,\zeta\)
\begin{description}
\item[(i)] \(|\uzi|< n_i< \omega\)
\item[(ii)] \(\uzi\subseteq\zeta\)
\item[(iii)] Let $B$ be a finite set of ordinals
and let $\zeta$ be such that $B\subseteq\zeta
<\lambda^+$.
Then $\{i:B\subseteq \uzi\}\in D$  
\item[(iv)] Coherency:
\(\gamma\in\uzi\Rightarrow\ugi=\uzi\cap\gamma\)

\end{description}  
\end{lemma}

Assuming the lemma, and
letting \(M=\{a_\zeta:\zeta<\lambda^+\}\) we now define, for each
\(\zeta\),
a function \(f_\zeta:\lambda \mapsto N\). 

Let $\Delta=\{\phi_{\alpha}:\alpha<\lambda\}$ and let 
$\{A_{\alpha}:\alpha <\lambda\}$ be a family witnessing the 
regularity of D. Thus for each $i$, the set $w_i=\{\alpha
:i\in A_{\alpha}\}$ is finite. Let $\Delta_i=\{\phi_{\alpha}
:\alpha\in w_i\}$, and let $u^{\zeta}_i, n_i$ be as in the
lemma.

We define a sequence of formulas essential to the proof:
suppose $\zeta<\lambda^+$ and $i<\lambda$. 
Let $m^{\zeta}_i=|\uzi|$ and let
\[\uzi=\{\xi_{\zeta,i,0},...,\xi_{\zeta,i,m^{\zeta}_i-1}\}\]
be an increasing enumeration of \(\uzi\).
(We adopt henceforth the convention that any 
enumeration of \(\uzi\) that is given
is an increasing enumeration.)
Let $\tzii$ be the $\Delta_i$-type of the tuple
$\langle 
a_{\xi_{\zeta,i,0}},...,
a_{\xi_{\zeta,i,m^{\zeta}_i-1}},a_\zeta\rangle$ in $M$. (So 
every $\phi(x_{\xi_{\zeta,i,0}},...,x_{\xi_{
\zeta,i,m^{\zeta}_i-1}},x_\zeta) 
\in \Delta _i $ or its negation occurs as a conjunct
of $\tzii$, according to whether 
$\phi(a_{\xi_{\zeta,i,0}},...,a_{\xi_{\zeta,i,m^{\zeta}_i-1}},a_\zeta)$
or 
$\neg\phi(a_{\xi_{\zeta,i,0}},...,a_{\xi_{\zeta,i,m^{\zeta}_i-1}},a_\zeta)$
holds in $M$.)
We define the formula $\tzi$ for each $i$ by downward
induction on $m^{\zeta}_i$ 
as follows: 
\smallskip

\noindent Case 1: $m^{\zeta}_i=n_i$. 
Let $\tzi=\bigwedge\tzii$. 
\smallskip

\noindent Case 2: $m^{\zeta}_i<n_i$. Let $\tzi$ be the 
conjunction of $\tzii$ and all formulas of the form
$\exists x_{m^{\epsilon}_i}\theta^{\epsilon}_i(x_0,...,x_
{m^{\zeta}_i-1},
x_{m^{\epsilon}_i})$, where $\epsilon$ satisfies
$u^{\epsilon}_i=u^{\zeta}_i\cup\{\zeta\}$ and hence
$m^{\epsilon}_i = m^{\zeta}_i+1$.

An easy induction shows that for a fixed $i<\lambda$,
the cardinality of the set
$\{\tzi:\zeta<\lambda^+\}$ is finite, using \(n_i\).

Let $i<\lambda$ be fixed. We define $f_\zeta(i)$ by induction
on $\zeta<\lambda^+$ in such a way that the following condition
remains valid:
\begin{enumerate}
\item[(IH)] If $\zeta^*<\zeta$ and $u^{\zeta^*}_i=
\{r_{\epsilon_1},...,r_{\epsilon_k}\}$, then 
$N\models\theta^{\zeta^*}_i
(f_{\epsilon_1}(i),...,f_{\epsilon_k}(i),f_{\zeta^*}(i)).$
\end{enumerate}

To define $f_\zeta(i)$, we consider different cases:
\smallskip

\noindent Case 1: $n_i=m^{\zeta}_i$.
\smallskip

\noindent Case 1.1: $n_i=0$.
Then $\tzi$ is the $\Delta_i$ type of the element $a_{\zeta}$.
But then 
\begin{eqnarray*}
M&\models&\tzi(a_{\zeta})\Rightarrow\\
M&\models&\exists x_0\tzi(x_0)\Rightarrow\\
N&\models&\exists x_0\tzi(x_0),\\
\end{eqnarray*}
where the last implication follows from the assumption that
$N$ satisfies the weakly $\Delta$-existential formulas 
holding in $M$.
Now choose an element $b\in N$ to witness this formula
and set $f_{\zeta}(i)=b.$

\noindent Case 1.2: $n_i>0.$
\smallskip
Let $\uzi=\{\xi_1,\ldots,\xi_{m^{\zeta}_i} \}$.
Since $m^{\zeta}_i=n_i$, 
the formula $\tzi$ is the $\Delta_i$-type of 
the elements $\{\xi_1,\ldots,\xi_{m^{\zeta}_i}\}$.
By assumption
$\gamma=\xi_{m^{\zeta}_i}$ is the maximum element of $\uzi$.  
Thus by coherency, $\ugi=\uzi\cap\gamma=
\{ \xi_1,\ldots,\xi_{m^{\zeta}_i-1}\}$. 
Since 
$\gamma<\zeta$, we know by the induction hypothesis that 
$$N\models\tgi
(f_{\xi_1}(i),\ldots,f_{\xi_{m^{\zeta}_i}}(i)).$$
By the formula construction $\tgi$ 
contains the formula $\exists x_{ m^{\zeta}_i}
\tzi(x_1,\ldots,x_{ m^{\zeta}_i})
$,
since $\uzi=\ugi\cup\{\gamma\}$ and since $m^{\gamma}_i<n_i.$
Thus $$N\models\exists x_{m^{\zeta}_i+1}
\tzi(f_{\xi_1}(i),\ldots,f_{\xi_{ m^{\zeta}_i}}(i),
x_{m^{\zeta}_i+1}).$$
As before choose an element $b\in N$ to witness this formula
and set $f_{\zeta}(i)=b.$

\noindent Case 2: $m^{\zeta}_i<n_i.$
Let $\uzi=\{\xi_1,\ldots,\xi_{m^{\zeta}_i} \}$.
We have that  $M\models\tzi(a_{\xi_1},\ldots,
a_{\xi_{m^{\zeta}_i}},a_{\zeta})$, and therefore  
$M\models \exists x_{m^{\zeta}_i +1}\tzi(a_{\xi_1},\ldots,
a_{\xi_{m^{\zeta}_i}},x_{m^{\zeta}_i+1})$. Let $\gamma=$ 
max($\uzi$)=$\xi_{m^{\zeta}_i}$.
By coherency, $\ugi=\uzi\cap\gamma$ and therefore since
$\gamma<\zeta$ by the induction hypothesis we have that
$$N\models\tgi
(f_{\xi_1}(i),\ldots,f_{\xi_{ m^{\zeta}_i -1}}(i),
f_{\gamma}(i)).$$
But then as in case 1.2 we can infer that
$$N\models\exists x_{m^{\zeta}_i+1}
\tzi(f_{\xi_1}(i),\ldots,f_{\xi_{m^{\zeta}_i }}(i),
x_{m^{\zeta}_i +1}).$$ As in case 1 choose an element $b\in N$ to 
witness this formula and set $f_{\zeta}(i)=b.$

It remains to be shown that the mapping
$a_{\zeta} \mapsto \langle f_{\zeta}(i):
i<\lambda \rangle /D$ satisfies the requirements of the
theorem, i.e. we must show, for all $\phi$ which is
in $\Delta$, or whose negation is in $\Delta$,
 $$M\models\phi(a_{\xi_1},\ldots,a_{\xi_k})
\Rightarrow \{i:N\models\phi(f_{\xi_1}(i),
\ldots, f_{\xi_k}(i))\}\in D.$$

\noindent So let such a $\phi$ be given,
 and 
%
%
suppose $M\models \phi (a_{\xi_1},
\ldots , a_{\xi_k})$.
Let $I_{\phi}=\{i:N\models\phi
(f_{\xi_1}(i),\ldots, f_{\xi_k}(i))
\}.$ We wish to show that $I_{\phi}\in D.$
Let $\alpha<\l$ so that $\phi$ is 
$\phi_\alpha$ or its negation.
It suffices to show that $A_{\alpha}\subseteq I_{\phi}.$ 
Let $\z<\l^+$ be such that $\{\xi_1,...,\xi_n\}\subseteq\z$.
By Lemma~\ref{main} condition (iii),
$\{i:\{\xi_1,...,\xi_n\}\subseteq \uzi\}\in D$.
So it suffices to show $$A_\alpha\cap
\{i:\{\xi_1,...,\xi_n\}\subseteq \uzi\}\subseteq I_\phi.$$
Let $i\in A_\alpha$ such that 
$\{\xi_1,...,\xi_n\}\subseteq \uzi$.
By the 
definition of $\theta^\z_i$ we know that
$N\models \theta^\z_i(f_{\xi_1}(i),\ldots, f_{\xi_k}(i)).$
But the 
$\Delta_i$-type of the tuple $\langle\xi_1,\ldots,\xi_k\rangle$
occurs as a conjunct of $\theta^\z_i$, and therefore
$N\models \phi(f_{\xi_1}(i),\ldots, f_{\xi_k}(i))$
\qed


\section{Proof of Lemma \ref{main}}


We now prove Lemma \ref{main}. We first prove a weaker
version in  which the filter is not given in advance:

\begin{lemma}\label{mainw}
Assume \(\gap\).
There exist sets 
\(\langle\uzi:\zeta<\lambda^+,i<cof(\lambda)\rangle\), integers $n_i$
and
 a regular filter $D$ on $\l$, generated by $\l$
sets,
such that
(i)-(iv) of Lemma~\ref{main} hold.
\end{lemma}

\proof
By \cite[Proposition 5.1, p. 149]{Sh269} the assumption \(\gap\) is
equivalent to:
\begin{itemize}
\item[\(\Box^{b^*}_\lambda:\)]
There is a $\l^+$-like linear order $L$, sets
$\langle\cza:a\in L,\z<cf(\l)\rangle$,
equivalence relations
$\langle\ez:\z<cf(\l)\rangle$, and
functions $\langle \fzab:\z<\l, a\in L, b\in L\rangle$
such that
\begin{enumerate}
\item[(i)] $\bigcup_\z\cza=\{b:b<_L a\}$ (an increasing
union in $\zeta$).
\item[(ii)] If \(b\in \cza\), then
$\czb=\{c\in\cza:c<_L b\}$.
\item[(iii)] $\ez$ is an equivalence relation on $L$
with $\le\l$ equivalence classes.
\item[(iv)] If $\z<\xi<cf(\l)$, then $\ex$ refines $\ez$.
\item[(v)] If $a\ez b$, then $\fzab$ is an order-preserving
one to one mapping from $\cza$ onto $\czb$ such that for 
$d\in \cza, d\ez \fzab(d)$.
\item[(vi)] If $\z<\xi<cf(\l)$ and $a\ex b$, then
$\fzab\subseteq\fxab$.
\item[(vii)] If $\fzab(a_1)=b_1$, then $f^\z_{a_1,b_1}
\subseteq\fzab$.
\item[(viii)] If $a\in \czb$ then $\neg \ez(a,b)$.

\end{enumerate}
\end{itemize}

This is not quite enough to prove Lemma \ref{mainw}, so we have
to work a little more. Let $$\Xi_{\z}=\{a/\ez:a\in L\}.$$
We assume, for simplicity, that $\z\ne\xi$ implies
$\Xi_\z\cap\Xi_\xi=\emptyset$. Define for $t_1,t_2\in\Xi_\z$:
$$t_1<_\z t_2\iff(\exists a_1\in t_1)(\exists a_2\in t_2)(a_1\in
C^\z_{a_2}).$$

\begin{proposition}\label{fact1}
$\langle \Xi_\z,<_\z\rangle$ is a tree order with \(cf(\lambda)\) 
as the set of levels.
\end{proposition}

\proof We need to show (a) $t_1<_{\zeta}t_2<_{\zeta}
t_3 $ implies $ t_1<_{\zeta}t_3$, and (b) 
$t_1<_{\zeta}t_3$ and $t_2<_{\zeta}t_3$ implies $t_1<
_{\zeta}t_2$ or $t_2<_{\zeta}t_1$ or \(t_1=t_2\). For the 
first, $t_1<_{\zeta}t_2$ implies there exists $a_1\in t_1$ 
and $a_2\in t_2$  
such that $a_1\in 
C^{\zeta}_{a_2}$. Similarly $t_2<_{\zeta}t_3$ implies
there exists $b_2\in t_2$ and $b_3\in t_3$  
such that $b_2 \in C^{\zeta}_{b_3}$. Now $a_2\ez b_2$ and 
hence we have the order preserving map $f^{\zeta}_{a_2,b_2}$
from $C^{\zeta}_{a_2}$ onto $C^{\zeta}_{b_2}$. Recalling
$a_1\in C^{\zeta}_{a_2}$, 
let $f^{\zeta}_{a_2,b_2}(a_1)=b_1$.
Then by (vi), $a_1\ez b_1$ and hence $b_1 \in t_1.$ But 
then $b_1 \in C^{\zeta}_{b_2}$ implies $b_1 \in C^{\zeta}_{b_3}$, 
by coherence and the fact that $b_2\in  C^{\zeta}_{b_3} .$ But then 
it follows that $t_1<_{\zeta}t_3$. 

Now assume
$t_1<_{\zeta}t_3$ and $t_2<_{\zeta}t_3$. Let 
$a_1\in t_1$ and $a_3\in t_3$ be such that $a_1\in C^{\zeta}_
{a_3}$, and similarly let $b_2$ and $b_3$ be such that 
$b_2 \in C^{\zeta}_{b_3}$. 
$a_3E^{\zeta}b_3$ implies we have the order preserving map
$f^{\zeta}_{a_3,b_3}$ from $C^{\zeta}_{a_3}$ to
$C^{\zeta}_{b_3}.$ Letting $f^{\zeta}_
{a_3,b_3}(a_1)=b_1$, we see that $b_1\in C^{\zeta}_{b_3}$.
If \(b_1<_L b_2\), then we have
$C^{\zeta}_{b_2}=C^{\zeta}_{b_3} \cap \{ c|c<b_2\} 
$ which implies $ b_1\in C^{\zeta}_{b_2}$, since, as $f^{\zeta}_
{a_3,b_3} $ is order preserving, $b_1<_Lb_2.$ Thus
$t_1<_{\zeta}t_2.$ The case $b_2<_L b_1$ is proved similarly,
and \(b_1=b_2\) is trivial.
\qed

For $a<_L b$ let
$$\xi(a,b)=\min\{\z:a\in \czb\}.$$
Denoting $\xi(a,b)$ by $\xi$, let
$$tp(a,b)=\langle a/E^{\xi},b/E^{\xi}\rangle.$$
If $a_1<_L...<_L a_n$, let
$$tp(\langle a_1,...,a_n\rangle)=
\{\langle l,m,tp(a_l,a_m)\rangle|1\le l<m\le n\}$$
and $$\Gamma=\{tp(\av):\av\in{}^{<\omega}L\}.$$
For $t=tp(\av), \av\in {}^nL$ we use $n_t$ to denote
the length of $\av$.

\begin{proposition}\label{fact2}
If $a_0<_L...<_L a_n$, then $$\max\{\xi(a_l,a_m):0\le l <m\le n\}=
\max\{\xi(a_l,a_n):0\le l< n\}.$$
\end{proposition}

\proof Clearly the right hand side is $\leq$ the left hand side.
To show the left hand side is $\leq$ the right hand side,
let $l<m <n$ be arbitrary. If $\xi(a_l,a_n)\le\xi(a_m,a_n)$,
then  $\xi(a_l,a_m)\le\xi(a_m,a_n)$. On the other hand,
if $\xi(a_l,a_n)>\xi(a_m,a_n)$,
then  $\xi(a_l,a_m)\le\xi(a_l,a_n)$. In either case
$\xi(a_l,a_m)\le\max\{\xi(a_k,a_n):0\le k< n\}$.
\qed

Let us denote max$\{\xi(a_l,a_n):0\le l< n\}$ by $\xi(\av)$.
We define on $\Gamma$ a two-place relation $\leq_{\Gamma}$
as follows:

$$ t_1 <_{\Gamma} t_2$$
if there exists
a tuple $\langle a_0,\ldots
a_{n_{t_2}-1} \rangle$ realizing $t_2$ such that some 
subsequence of the tuple realizes $t_1$. 

Clearly,
$\langle\Gamma,\leq_{\Gamma}\rangle$ is a directed partial order.

\begin{proposition}\label{fact5} For $t\in \Gamma$,
$t=tp(b_0,\ldots b_{n-1})$ and
$a\in L$, there exists at most one $k <n$ such that
$b_{k}E^{\xi(b_0,\ldots,b_{n-1})}a$. 
\end{proposition}

\proof Let $\z=\xi(b_0,\ldots,b_{n-1})$ and let 
$b_{k_1}\ne b_{k_2}$ be such that $b_{k_1}E^{\zeta}
a$ and $b_{k_2}E^{\zeta}a$, $k_1,k_2\le n-1.$ Without
loss of generality, assume $b_{k_1}<b_{k_2}.$ 
Since $E^{\zeta}$ is an equivalence relation,
$b_{k_2}E^{\zeta}b_{k_1}$ and thus we have an 
order preserving map $f^{\z}_{b_{k_2},b_{k_1}}$ from
$C^{\z}_{b_{k_2}}$ to $C^{\z}_{b_{k_1}}$. Also 
$b_{k_1}\in C^{\z}_{b_{k_2}}$, by the definition
of $\z$ and by coherence, and therefore 
$f^{\z}_{b_{k_2},b_{k_1}}(b_{k_1})E^{\z}b_{k_1}.$ 
But this contradicts (viii), since  
$f^{\z}_{b_{k_2},b_{k_1}}(b_{k_1})\in C^{\z}_{b_{k_1}}$.
\qed

\begin{definition}\label{ww} For $t\in \Gamma$,
$t=tp(b_0,\ldots b_{n-1})$ and
$a\in L$ suppose there exists $k <n$ such that
$b_{k}E^{\xi(b_0,\ldots,b_{n-1})}a$. Then let
\(u^a_t = \{f^{\z (b_0,\ldots ,b_{n-1})}_{a,b_{k}}(b_l):l<k\} \)
Otherwise, let \(u^a_t=\emptyset\).
\end{definition}

Finally, let $D$ be the filter on $\Gamma$ generated
by the $\l$ sets 
$$\Gamma_{\ge t^*}=\{t\in \Gamma:t^*<_L t\}.$$
We can now see that the sets $\uat$, the numbers
$n_t$ and the filter $D$ satisfy conditions
(i)-(iv) of Lemma~\ref{main} with $L$ instead of
\(\lambda^+\):
Conditions (i) and (ii) are trivial in this case.
Condition (iii) is verified as follows: Suppose $B$
is finite. Let $a\in L$
be such that $(\forall x\in B)(x<_L a)$.
Let $\av$ enumerate $B\cup\{a\}$ in increasing order
and let $t^*=tp(\av)$. Clearly
$$t\in\Gamma_{\ge t^*}\Rightarrow B\subseteq \uat.$$
Condition (iv) follows directly from Definition~\ref{ww}
and Proposition~\ref{fact5}.

To get the Lemma on $\l^+$ we observe that since $L$ is
$\l^+$-like, we can assume that $\langle\l^+,<\rangle$
is a submodel of $\langle L, <_L\rangle$. Then we define
$v^\alpha_t=u^\alpha_t\cap\{\beta:\beta<\alpha\}$.
Conditions (i)-(iv) of Lemma~\ref{mainw} are still satisfied.
Also having \(D\) a filter of \(\Gamma\) instead of
\(\lambda\) is immaterial as \(|\Gamma|=\lambda\). \qed

Now back to the proof of Lemma~\ref{main}. Suppose
$\uzi,n_i$ and $D$ are as in Lemma~\ref{mainw}, and
suppose $D'$ is an arbitrary regular filter on $\l$.
Let $\{A_{\alpha}:\alpha <\l \}$ be a family of sets witnessing
the regularity of $D'$, and let $\{Z_{\alpha}:\alpha <\l \}$
be the family generating $D.$ We define a function $h:\l
\rightarrow \l$ as follows. Suppose $i < \l.$ Then
let $$h(i)\in \bigcap\{Z_{\alpha}|i\in A_{\alpha}\}.$$
Now define $v^{\z}_{\alpha}=u^{\z}_{h(\alpha)}$. Define also
 $n_{\alpha}=n_{h(\alpha)}.$ Now the sets $v^{\z}_{\alpha}$
and the numbers 
$n_{\alpha}$ satisfy the conditions of Lemma~\ref{main}.
\qed


\section{Is $\Box^{b^*}_\l$  needed for Lemma \ref{mainw}?}


In this section we 
show that the conclusion of Lemma~\ref{mainw}
(and hence of Lemma~\ref{main}) implies $\Box^{b^*}_\l$
for singular strong limit \(\lambda\). By
\cite[Theorem 2.3 and Remark 2.5]{Sh269}, $\Box^{b^*}_\l$
is equivalent, for singular strong limit \(\lambda\),
to the following principle:

\def\S{{\cal S}}

\begin{itemize}
\item[\(\S_\lambda:\)]
There  are sets
$\langle C^i_a:a< \lambda^+,i<cf(\l)\rangle$
such that
\begin{enumerate}
\item[(i)] If $i<j$, then $C_a^i\subseteq C_a^j$.
\item[(ii)] $\bigcup_i C^i_a=a$.
\item[(iii)] If \(b\in C^i_a\), then
$C^i_b=C^i_a\cap b$.
\item[(iv)] $\sup\{otp(C^i_a):a<\l^+\}<\l$.
\end{enumerate}
\end{itemize}

Thus it suffices to prove:

\begin{proposition}
Suppose the sets $\uzi$ and
the filter $D$ are as given by Lemma~\ref{mainw}
and \(\lambda\) is a limit cardinal.
Then $\S_\lambda$ holds.
\end{proposition}

\proof
Suppose ${\cal A}=\{A_{\alpha}:\alpha <\l \}$ is a family of sets 
generating $D$. W.l.o.g.,  ${\cal A}$ is closed under finite intersections.
Let \(\l\) be the union of the
increasing sequence \(\langle \l_\alpha :\alpha<cf(\l)\rangle\), where \(\l_0\ge\omega\).
Let the sequence $\langle \Gamma_\alpha:\alpha<cf(\lambda)\rangle$
satisfy: 
\begin{enumerate}
\item[(a)] $|\Gamma_\alpha| \le
\l_\alpha$
\item[(b)] $\Gamma_\alpha$ is continuously
increasing in $\alpha$ with $\l$ as union
\item[(c)] If $\beta_1,...,\beta_n\in\Gamma_\alpha$,
then there is $\gamma\in\Gamma_\alpha$ such that
$$A_\gamma=A_{\beta_1}\cap...\cap A_{\beta_n}.$$
\end{enumerate}
The sequence $\langle \Gamma_\alpha:\alpha<cf(\lambda)\rangle$ enables
us to define a sequence that will witness $\S_\l$.
For $\alpha<cf(\l)$ and $\z<\l^+$, let
\def\vza{V^\alpha_\z}
\def\vxa{V^\alpha_\xi}
$$\vza=\{\xi<\z:(\exists\gamma\in\Gamma_\alpha)
(A_\gamma\subseteq\{i:\xi\in \uzi\})\}.$$

\begin{lemma}\label{box}
\begin{enumerate}
\item[(1)] $\langle \vza:\alpha<\l\rangle$ is a continuously
 increasing sequence of subsets of $\zeta$, $|\vza|\le \l_\alpha$,
and $\bigcup\{\vza:\alpha<cf(\l)\}=\zeta$.
\item[(2)] If $\xi\in\vza$, then $V^\alpha_\xi=\vza\cap\xi$.
\end{enumerate}
\end{lemma}

\proof (1) is a direct consequence of the definitions.
(2) follows from the respective property of the sets
\(\uzi\). \qed

\begin{lemma}\label{bbox}
$\sup\{otp(\vza):\z<\l^+\}\le \l_{\alpha}^+$.
\end{lemma}

\proof By the previous Lemma,  $|\vza|\le\l_\alpha$.
Therefore $otp(\vza)<\l_\alpha^+$
and the claim follows. \qed

The proof of the proposition is complete: (i)-(iii) follows from
Lemma~\ref{box}, (iv) follows from Lemma~\ref{bbox} and the assumption
that \(\l\) is a limit cardinal. \qed

More equivalent conditions for the case $\l$ singular strong limit, $D$
a regular ultrafilter on $\l$, are under preparation.



\def\gg{\gamma}
\newcommand\EFG[1]{{\rm EFG}_{#1}}
\def\I{{\rm I}}
\def\II{{\rm II}}
\def\la{\langle}
\def\ra{\rangle}\def\ga{\alpha}
\def\gb{\beta}\def\se{\subseteq}
\def\uai{u^{\alpha}_i}
\def\uzi{u^{\zeta}_i}
\def\uei{u^{\epsilon}_i}
\def\ugi{u^{\gamma}_i}
\def\ez{E^{\zeta}}
\def\ex{E^\xi}
\def\cza{C^{\zeta}_a}
\def\czb{C^{\zeta}_b}
\def\z{\zeta}
\def\l{\lambda}
\def\fzab{f^\z_{a,b}}
\def\fxab{f^\xi_{a,b}}
\def\av{\vec{a}}
\def\bv{\vec{b}}
\def\max{{\rm max}}
\def\min{{\rm min}}
\def\uat{u^a_t}
\def\utb{u^b_t}
\def\vat{v^a_t}
\def\vat{v^a_t}\def\qed{$\Box$\medskip}
\def\tzi{\theta^{\zeta}_i}
\def\tei{\theta^{\epsilon}_i}
\def\tzii{\bar{\theta}^{\zeta}_i}
\def\tgi{\theta^{\gamma}_i}
\def\proof{\noindent{\bf Proof.}\hspace{2mm}}

\section{Ehrenfeucht-Fra\"\i ss\'e-games}

Let \(M\) and \(N\) be two first order structures of the
same vocabulary \(L\). 
All vocabularies are assumed to be relational. 
The  
{\em Ehrenfeucht-Fra\"\i ss\'e-game of length \(\gg\) 
of \(M\) and \(N\)}
denoted  by \(\EFG{\gg}\) is defined as follows: There are
 two players
called \(\I\) and \(\II\). First \(\I\) plays \(x_0\) and 
then \(\II\) plays
\(y_0\). After this
 \(\I\) plays \(x_1\), and \(\II\) plays \(y_1\), and
so on. If  \(\la(x_{\gb},y_{\gb}):\gb<\ga\ra\) has been 
played
and \(\ga<\gg\), then  \(\I\) plays \(x_{\ga}\) after 
which \(\II\)
plays \(y_{\ga}\). Eventually a sequence
\(\la(x_{\gb},y_{\gb}):\gb<\gg\ra\) has been played. The rules
of the game say that both players have to play elements
of \(M\cup N\). Moreover, if 
\(\I\) plays his \(x_{\gb}\) in \(M\) (\(N\)), then \(\II\) 
has to play
his \(y_{\gb}\) in \(N\) (\(M\)).
Thus the sequence \(\la(x_{\gb},y_{\gb}):\gb<\gg\ra\)
determines a relation \(\pi\se M\times N\). 
Player \(\II\) wins this round of the game if \(\pi\) is
a partial isomorphism. Otherwise \(\I\) wins. 
The notion of winning strategy is defined in the usual manner.
We say that a player {\em wins} \(\EFG{\gg}\) if he has 
a winning
strategy in \(\EFG{\gg}\).

Note that if II has a winning strategy in \(\EFG{\gg}\)
on \(M\) and \(N\), where \(M\) and \(N\) are of size
\(\le|\gg|\), then \(M\cong N\).

Assume $L$ is of cardinality $\le\l$
and for each $i<\l$ let $M_i$ and $N_i$  are 
elementarily equivalent $L$-structures. 
Shelah proved in \cite{class}
that if  $D$ is a regular filter on $\l$,
then Player II has a winning strategy in the game
$\EFG{\gamma}$ on $\prod_iM_i/D$
and $\prod_iN_i/D$ for each \(\gamma<\lambda^+\).
 We show that
under a stronger assumption, II has a winning
strategy even in the game $\EFG{\lambda^+}$.
This makes a big difference because, assuming
the models \(M_i\) and \(N_i\) are
of size \(\le\lambda^+\), \(2^{\lambda}=\lambda^+\), and the models
$\prod_iM_i/D$
and $\prod_iN_i/D$ are of size \(\le\lambda^+\). Then by the
remark above, if II has a winning strategy in $\EFG{\lambda^+}$,
the reduced powers are actually isomorphic.
Hyttinen \cite{hytt} proved this under the assumption
that the filter is, in his terminology, semigood.

\begin{theorem} \label{second} Assume $\gap$.
Let $L$ be a vocabulary of cardinality $\le\l$
and for each $i<\l$ let $M_i$ and $N_i$ be
two elementarily equivalent $L$-structures. If $D$ is a regular filter on $\l$,
then Player II has a winning strategy in the game
$\EFG{\l^+}$ on $\prod_iM_i/D$
and $\prod_iN_i/D$.
\end{theorem}

\proof
We use Lemma~\ref{main}.
If $i<\l$, then, since  $M_i$ and $N_i$ are 
elementarily equivalent, Player II has a
winning strategy $\sigma_i$ in the game $\EFG{n_i}$
on $M_i$ and $N_i$.
We will use the set $\uzi$ to put these short
winning strategies together into one long
winning strategy.

A ``good'' position is a sequence 
$\langle (f_\z,g_\z):\z<\xi\rangle$, where $\xi<\l^+$,
and for all $\z<\xi$ we have
$f_\z\in\prod_iM_i$, $g_\z\in\prod_iN_i$, and
if $i<\l$, then 
$\langle (f_\epsilon(i),g_\epsilon(i)):\epsilon\in\uzi\cup\{\z\}\rangle$
is a play according to $\sigma_i$.

Note that in a good position the equivalence classes
of the functions \(f_\zeta\) and \(g_\zeta\) determine
a partial isomorphism of the reduced products.
The strategy of player II is to keep the position of the game ``good'',
and thereby win the game.
Suppose \(\xi\) rounds have been played and II has been able
to keep the position ``good''. Then player I plays \(f_\xi\).
We show that player II can play \(g_\xi\) so that
$\langle (f_\z,g_\z):\z\le\xi\rangle$ remains ``good''.
Let \(i<\lambda\). Let us look at
$\langle (f_\epsilon(i),g_\epsilon(i)):\epsilon\in
u^{\xi}_i\rangle$.
We know that this is a play according to the strategy \(\sigma_i\)
and \(|u^{\xi}_i|<n_i\). Thus we can play one more move
in \(EF_{n_i}\) on \(M_i\) and \(N_i\) with player I playing
\(f_\xi(i)\). Let \(g_\xi(i)\) be the answer of II in
this game according to \(\sigma_i\). The values
\(g_\xi(i)\), \(i<\lambda\), constitute the function
\(g_\xi\). We have showed that II can  maintain a ``good''
position.
\qed

\begin{corollary} Assume GCH and $\l$ regular (or just $\gap$
and $2^\l=\l^+$). 
Let $L$ be a vocabulary of cardinality $\le\l$
and for each $i<\l$ let $M_i$ and $N_i$ be
two elementarily equivalent
$L$-structures. If $D$ is a regular filter on $\l$,
then  $\prod_iM_i/D\cong\prod_iN_i/D$.
\end{corollary}

\end{document}